# Non-singularity of the generalized logit dynamic with an application to fishing tourism

Hidekazu Yoshioka[1, a)]

[1]*Japan Advanced Institute of Science and Technology, 1-1 Asahidai, Nomi, Ishikawa 923-1292, Japan*

[a)] *Corresponding author: yoshih@jaist.ac.jp*

**Abstract.** Generalized logit dynamic defines a time-dependent integro-differential equation with which a Nash equilibrium of an iterative game in a bounded and continuous action space is expected to be approximated. We show that the use of the exponential logit function is essential for the approximability, which will not be necessarily satisfied with functions with convex exponential-like functions such as *q*-exponential ones. We computationally analyze this issue and discuss influences of the choice of the logit function through an application to a fishing tourism problem in Japan.

## INTRODUCTION

A Nash equilibrium is a strategy profile of an iterative game such that no player can do better by deviating, which has been pointed out to be not robust against perturbation[1]. Logit dynamics have been utilized for efficiently approximating Nash equilibria in the framework of evolutionally systems, with which optimal actions of agents are obtained as smooth, non-singular functions but not as a set of Dirac deltas[2]. This smoothing property of logit dynamics emerges as an equilibrium being robust against perturbation, often called the logit equilibrium[2].

We formally introduce a logit dynamic following the literature[3]. The time $t$ is a non-negative parameter. For a compact space $\Omega$, such as a bounded closed 1-D interval, parameterized by $x$, denote a payoff function of an agent choosing the action $x$ by $\varphi(x;\mu)$ that possibly depends on (conditional) expectations with respect to $\mu$. A time-dependent probability measure $\mu = \mu_t(\mathrm{d}x)$ on $\Omega$ represents the distribution of actions of agents parameterized by $x$. Set a smoothing parameter $\eta > 0$ and an initial distribution $\mu_0(\mathrm{d}x)$, the logit dynamic that governs the temporal evolution of $\mu$ is the integro-differential equation

$$\frac{\mathrm{d}\mu}{\mathrm{d}t}(A) = L_\eta \mu(A) - \mu(A), \; t > 0 \tag{1}$$

for any Borel measurable set $A$ of $\Omega$, where

$$L_\eta \mu(A) = \int_A \exp\left(\frac{\varphi(y;\mu)}{\eta}\right)\mathrm{d}y \bigg/ \int_\Omega \exp\left(\frac{\varphi(y;\mu)}{\eta}\right)\mathrm{d}y. \tag{2}$$

It has been proven that, under a mild condition, the integro-differential equation (1) admits a unique measure-valued solution, which eventually satisfies the equality

$$L_\eta \mu(A) = \mu(A) \text{ as } t \to +\infty \tag{3}$$

for any Borel measurable set $A$ of $\Omega$ [3]. It was also proven that the steady solution $\mu$ to (3) converges to a Nash equilibrium with the payoff function $\varphi(x;\mu)$ under $\eta \to +0$.

Logit dynamics have been utilized for approximating stable Nash equilibria, while the dynamics themselves also are interesting because they represent iterative games under action uncertainty of agents. Indeed, the term $L_\eta \mu$ is interpreted as the maximizer of an optimization problem:

$$L_\eta \mu = \underset{\phi \geq 0, \int_\Omega \phi(x)\mathrm{d}x=1}{\arg\max} \left\{ \underbrace{\int_\Omega \phi(x)\varphi(x;\mu)\mathrm{d}x}_{\text{Objective}} - \eta \underbrace{\int_\Omega \phi(x)\ln\phi(x)\mathrm{d}x}_{\text{Relative entropy}} \right\}, \quad (4)$$

which is a maximization problem of an objective (first term) subject to the model uncertainty through a Lagrangian multiplier $\eta$ (second term). Notice that $\mu$ is a given probability measure in (4). The smoothing parameter $\eta$ serves as a width (i.e., uncertainty) of the profile of $\mu$, thus avoiding the singularity corresponding to pure strategies. It has been pointed out that the term $L_\eta \mu$ needs not arise from (4) but can be formulated through a variant where the relative entropy is replaced by some other statistical divergence such as the Tsallis or Rényi one[4]. Such a logit dynamic is called a generalized logit dynamic. It was pointed out in the literature[4] that possibly different Nash equilibria emerge with different divergences; however, there was no clear explanation why this phenomenon occurs.

The objective of this work is to present a generalized logit dynamic motivated from a fishing tourism problem in a mountainous region in Japan. We theoretically show that the problem has a pure Nash equilibrium. Then, we demonstrate that the associated generalized logit dynamics, in the context of a Tsallis divergence generalizing (4)[5], with different shape parameter values lead to different steady $\mu$. We further show that the pure Nash equilibrium cannot be approximated by the Tsallis divergence unless it exactly coincides with the classical relative entropy of (4). We also discuss fisheries implications of this phenomenon.

## MAIN RESULTS

### Generalized Logit Dynamic

Hereafter, we choose $\Omega = [0,1]$ for sake of simplicity, but any other bounded 1-D intervals can be chosen if preferred. Introduce a shape parameter $q \geq 1$ and consider an extended version of (4) in the sense of Tsallis:

$$L_{\eta,q} \mu = \underset{\phi \geq 0, \int_\Omega \phi(x)\mathrm{d}x=1}{\arg\max} \left\{ \underbrace{\int_\Omega \{\phi(x)\}^q \varphi(x;\mu)\mathrm{d}x}_{\text{Objective}} - \eta \underbrace{\int_\Omega \frac{1}{1-q}\left(1-q+q\phi(x)-\{\phi(x)\}^q\right)\mathrm{d}x}_{\text{Tsallis divergence}} \right\}, \quad (5)$$

which formally reduces to (4) as $q \to 1$ from above. Assume that $\varphi(x;\mu)$ is not positive. The maximizer of (5), which is $L_{\eta,q}$ in the left-hand side of (5), is explicitly obtained as

$$L_\eta \mu(A) = \int_A \exp_q\left(\frac{\varphi(y;\mu)}{\eta}\right) \mathrm{d}y \Big/ \int_\Omega \exp_q\left(\frac{\varphi(y;\mu)}{\eta}\right) \mathrm{d}y \quad (6)$$

for any Borel measurable set $A$ of $\Omega$ [3], where $\exp_q$ is the $q$-exponential function:

$$\exp_q(z) = \left(1+(1-q)z\right)^{\frac{-1}{q-1}}, \quad z \leq 0. \quad (7)$$

For any Borel measurable set $A$ of $\Omega$, the generalized logit dynamic subject to an initial condition $\mu_0$ is set as

$$\frac{\mathrm{d}\mu}{\mathrm{d}t}(A) = L_{\eta,q}\mu(A) - \mu(A), \quad t > 0 \quad (8)$$

and its steady version as

$$L_{\eta,q}\mu(A) = \mu(A). \quad (9)$$

## Case Study

We set the payoff and explain its motivation. Here, $x$ is the frequency of fishing activities. For any $x, y \in \Omega$, set

$$f(x, y) = \underbrace{-\left(\frac{1}{x(1-x)}\right)^\gamma}_{\text{Penalizing too sparse/dense fishing activities}} \underbrace{-\alpha \exp(-\beta(x-y))}_{\text{Penalizing lower fishing utility of the agent } x \text{ than the other agents } y} \tag{10}$$

with parameters $\alpha, \beta, \gamma > 0$, and then

$$\varphi(x; \mu) = \int_\Omega f(x, y) \mu(\mathrm{d}y). \tag{11}$$

This payoff $\varphi$ is motivated from our recent field surveys during the spring in 2023 around Shiramine Village in Hakusan mountainous region (hereafter simply called Hakusan), Ishikawa Prefecture, Japan. The Tedori River system in Hakusan is famous with mountain stream fishing of Japanese char, commonly called Iwana in Japan, during spring to summer in each year. Hakusan-Shiramine Fishery Cooperatives in Shiramine Village authorize inland fisheries in Hakusan. The upstream reaches of the Tedori River and their watershed area, near Shiramine Village, are characterized by steep slopes at a danger of collapse (called Kuzure). The union president of Hakusan-Shiramine Fishery Cooperatives claimed that tourists visiting this region for stream fishing can find dangerous areas at a danger of collapse earlier than residents, and their early warning is helpful for disaster prevention of residents of the village (personal communication on May 31, 2023 at the resultant Koyo-mon in Shiramine Village). By contrast, overfishing should be prevented to sustain the population of Iwana in the Tedori River system.

Now, we can explain each term of (10). The first term of (10) represents the disutility caused by too sparse or dense fishing activities. The second term of (10) represents the minus of a relative fishing utility.

## Numerical and Mathematical Analyses

This model is very simple, while it is sufficient to tackle the convergence issue explained in **INTRODUCTION**. A straightforward calculation shows that the pure Nash equilibrium of the iterative game maximizing the payoff in (11) is given by the unique positive solution $x = \hat{x} \in (0,1)$ to ($\delta_{\{\cdot = x\}}$ is the Dirac delta concentrated at $x$)

$$0 = -\frac{\mathrm{d}}{\mathrm{d}x}\left(\frac{1}{x(1-x)}\right)^\gamma + \alpha\beta \exp(-\beta x) \int_0^1 \exp(\beta y) \delta_{\{y=x\}} \mathrm{d}y \quad \text{or equivalently} \quad \frac{\mathrm{d}}{\mathrm{d}x}\left(\frac{1}{x(1-x)}\right)^\gamma = \alpha\beta, \tag{12}$$

which can be computed by a common fixed-point method.

**Figure 1** in the next page shows the computed numerical solutions to (9) where the domain $\Omega$ is uniformly discretized by 200 cells on which the integro-differential equation (8) is temporally integrated by a Forward Euler method with the time increment of 0.1. We set $(\alpha, \beta, \gamma) = (1,1,1)$. **Figures 1(a)-(b)** demonstrate that $\mu$ concentrates on the pure Nash equilibrium $x = \hat{x}$ as $\eta \to +0$ if $q = 1$, while it is not if $q = 1.2$ (the same seems apply to the other $q > 1$); namely, the classical relative entropy is essential for the approximability viewpoint. Note that if $q = 1.2$ the mode of the limit distribution captures the pure Nash strategy according to **Figure 1(b)**.

We explain the convergence phenomena in **Figures 1(a)-(b)** theoretically. Assume that there exists a probability density function $p = p(x)$ of $\mu$ that is absolutely continuous with respect to the Lebesgue measure $\mathrm{d}x$. Notice that

$$\int_0^1 \exp(\beta y) \mu(\mathrm{d}y) = \int_0^1 \exp(\beta y) p(y) \mathrm{d}y \in [1, \exp(\beta)] \text{ is a constant}. \tag{13}$$

Then, for sufficiently small $\eta > 0$, with $q > 1$, we have (with some $\mu$)

$$\exp_q\left(\frac{\varphi(x;\mu)}{\eta}\right) = \left(1 + (1-q)\frac{\varphi(x;\mu)}{\eta}\right)^{\frac{-1}{q-1}} \sim \left(\frac{-\varphi(x;\mu)}{\eta}\right)^{\frac{-1}{q-1}}, \tag{14}$$

from which we obtain (notice the division by $\eta^{-1} > 0$ in the right-hand side of the equation just below)

$$p(x) \sim \left(\frac{-\varphi(x;\mu)}{\eta}\right)^{\frac{-1}{q-1}} \Big/ \int_\Omega \left(\frac{-\varphi(y;\mu)}{\eta}\right)^{\frac{-1}{q-1}} \mathrm{d}y = (-\varphi(x;\mu))^{\frac{-1}{q-1}} \Big/ \int_\Omega (-\varphi(y;\mu))^{\frac{-1}{q-1}} \mathrm{d}y. \tag{15}$$

By contrast, with $q = 1$, due to a classical saddle-point argument we obtain

$$p(x) = \exp\left(\frac{\varphi(x;\mu)}{\eta}\right) \Big/ \int_\Omega \exp\left(\frac{\varphi(y;\mu)}{\eta}\right) \mathrm{d}y \to \delta_{\{x=\hat{x}\}} \text{ in the sense of distributions.} \tag{16}$$

The mathematical analysis presented above shows that the difference of the growth speed between the exponential and $q$-exponential functions triggered the difference between the steady solutions.

In a fisheries viewpoint, using the Tsallis divergence implies that the individual anglers assume a strong diversity about the fishing (dis-)utilities among them, such that there exists no pure Nash strategy.

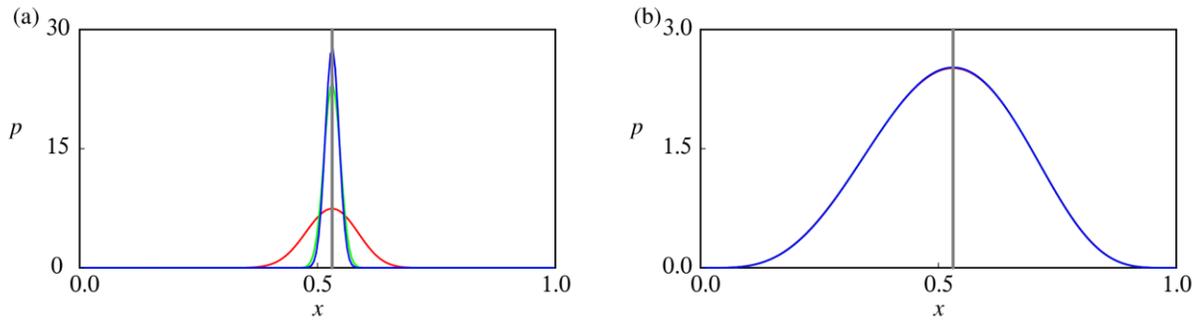

**FIGURE 1.** The computed equilibria solving (9) for (a) $q = 1$ and (b) $q = 1.2$ with different values of $\eta$: in panel (a), 0.1 (red), 0.01 (green), 0.007 (blue), and in panel (b), 0.01 (red), 0.001 (green), 0.0001 (blue). The grey line represents the pure Nash equilibrium $\hat{x} = 0.531$. The initial condition is the uniform distribution on $\Omega$.

## FUTURE PERSPECTIVES

We focused on the single-population problem in this work, while the proposed formulation can be extended to multi-population ones[5] in a straightforward manner. The payoff $\varphi$ can also be generalized. In the real world, fishing activities in a mountainous region have seasonality. Considering time-dependent iterative game problems are currently under investigations by the author and his coworkers. Another issue will be to analyze whether the mode of the Tsallis case equals the pure Nash equilibrium or not.

## ACKNOWLEDGMENTS

This research was supported by the Japan Society for the Promotion of Science, grant numbers 22K14441 and 22H02456. We thank the union president of Hakusan-Shiramine Fishery Cooperatives for motivating this work and providing useful comments on inland fisheries in the Tedori River system.